\documentclass{amsart}
\usepackage{amsmath,amssymb,graphicx,setspace,verbatim}
\usepackage{color, tikz, url}
\input xy
\xyoption{all}
\newtheorem{prop}{Proposition}[section]

\newtheorem{thm}[prop]{Theorem}
\newtheorem{lemma}[prop]{Lemma}

\newcommand{\Aut}{\mathrm{Aut}}

\newcommand{\po}{\mathcal{P}}
\newcommand{\qo}{\mathcal{Q}}
\newcommand{\fl}{\mathcal{F}}

\begin{document}
\title[A new algorithm for chiral polytopes]{A new algorithm to classify chiral polytopes with a given automorphism group}

\author[F. Buekenhout]{Francis Buekenhout}
\address{Francis Buekenhout, Universit\'e Libre de Bruxelles, D\'epartement de Math\'ematique, C.P.216 - Alg\`ebre et Combinatoire, Boulevard du Triomphe, 1050 Brussels, Belgium
}
\email{Francis.Buekenhout@ulb.ac.be}

\author[D. Leemans]{Dimitri Leemans}
\address{Dimitri Leemans, Universit\'e Libre de Bruxelles, D\'epartement de Math\'ematique, C.P.216 - Alg\`ebre et Combinatoire, Boulevard du Triomphe, 1050 Brussels, Belgium
}
\email{dleemans@ulb.ac.be}

\author[P. Tranchida]{Philippe Tranchida}
\address{Philippe Tranchida, Department of Mathematical Sciences, KAIST, 291 Daehak-ro Yuseong-gu, Daejeon, 34141, South Korea}
\email{ptranchi@kaist.ac.kr}

\date{}
\maketitle
\begin{abstract}
We present a new algorithm to compute all the chiral polytopes that have a given group $G$ as full automorphism group. This algorithm uses a new set of generators that characterize the group, all of them except one being involutions. It permits to compute all chiral polytopes of groups that were previously unreachable by other known algorithms.
\end{abstract}
\section{Introduction}
Abstract polytopes are combinatorial structures that generalize the face lattice of convex polytopes. 
Abstract polytopes are however not necessarily convex, nor necessarily realizable.
An abstract polytope is a partially ordered set $(\mathcal P,\leq)$, whose elements are called faces, and that satisfies some extra axioms (see Section~\ref{basics} for the definitions).

Every rank two abstract polytope is isomorphic to a regular polygon.
Abstract polytopes of rank 3 are essentially non-degenerate maps (i.e. 2-embeddings of maps on surfaces, see for instance~\cite[§3.2]{CM1972} for a precise definition): while every rank 3 polytope is indeed a map, the converse is not true as a map can fail the so-called diamond condition. In that case we call such a map degenerate.

Among abstract polytopes, those that possess a high degree of symmetry have been the most studied.
The ones having maximal possible rotational and reflectional symmetries are called regular while those having maximal rotational symmetries but no reflections are called chiral.

There are three atlases of regular polytopes; all of them classify the polytopes by their automorphism groups. The first of these atlases, \cite{mihAtlas}, contains information about all regular polytopes with automorphism group of size $n$, where $n$ is at most 2000, and not equal to 1024 or 1536. It has been extended by Marston Conder\footnote{See \url{https://www.math.auckland.ac.nz/~conder/}}. The second atlas, \cite{LV2005}, contains all regular polytopes whose automorphism group is an almost simple group $\Gamma$ such that $S \leq \Gamma \leq \Aut(S)$, where $S$ is a simple group and $\Gamma$ is a group of order less than 1 million appearing in the Atlas of Finite Groups (\cite{Con85}). The third atlas \cite{mihotherAtlas} extends the second atlas to sporadic groups. It has also been extended by Leemans and Mark Mixer in~\cite{Algo}. These atlases have given rise to several conjectures that have been proven over the years. We refer to~\cite{survey} for a survey of the results obtained on almost simple groups.

In~\cite{HHL}, Hartley, Hubard and Leemans produced two atlases of chiral polytopes by designing new algorithms to classify abstract chiral polytopes. One of them is given a group as input and produces all abstract chiral polytopes, up to isomorphism, that have that group as full automorphism group.
That algorithm was based on the description of the automorphism group of an abstract chiral polytope given by Schulte and Weiss in~\cite{Chiral}.
The recent paper of Fernandes, Leemans and Weiss on hypertopes~\cite{Hypertopes} (that are actually thin residually connected incidence geometries) suggests a better algorithm, based essentially on searching for tuples of involutions together with a non-involutory element of the group, instead of tuples of non-involutory elements of the group.
We describe in this paper a new algorithm following the ideas of~\cite{Hypertopes} that greatly outperforms the algorithm of Hartley, Hubard and Leemans.
The paper is organised as follows.
In Section~\ref{basics}, we give the necessary definitions and notations needed to understand this paper.
In Section~\ref{groups}, we recall how to construct chiral polytopes from groups.
In Section~\ref{algo}, we explain our new algorithm and show that it produces all the chiral polytopes on which a given group $G$ acts as automorphism group. We also present an implementation of this algorithm, written in {\sc Magma}\cite{Magma}. Finally we give some computing times showing that our new algorithm is much faster than the Hartley-Hubard-Leemans algorithm and that it permits to collect experimental data for much larger groups than the previous algorithm.

\section{Basic definitions and notations}\label{basics}
We start by recalling the basic theory of abstract polytopes and regular polytopes, as it was done in~\cite{HHL}. For details, we refer the reader to \cite{ARP}.

An {\em (abstract) polytope} of rank $n$ or an {\em  $n$-polytope}  is a partially ordered set (or poset) ${\mathcal P}$ endowed with a strictly monotone rank function having range $\{-1, \ldots, n\}$. 
The elements of $\po$ are called {\em faces}. For $\,0 \le j < n$, a face of rank $j$ is called a {\em $j$-face} and the faces of rank $\,0,1$ and $n-1$ are called the {\em vertices}, {\em edges} and {\em facets} of the polytope, respectively. 
We ask that ${\mathcal P}$ has a smallest face $F_{-1}$, and a greatest face $F_n$ (called the {\em improper faces} of $\po$), and that each {\em flag} (that is, each maximal chain of the order, also called {\em chamber} in incidence geometry) of ${\mathcal P}$ contain exactly $n+2$ faces. 
We denote by ${\mathcal F(\mathcal P)}$ the set of all flags of ${\mathcal P}$.
Two flags are said to be  {\em adjacent} if they differ by exactly one face,  they are $j$-adjacent, if the rank of the face they differ on is precisely $j$.
We also require ${\mathcal P}$ to be {\em strongly flag-connected}, that is, any two flags $\Phi, \Psi \in \fl(\po)$  can be joined by a sequence of flags $\Phi = \Phi_0, \Phi_1, \ldots, \Phi_k=\Psi$ such that 
each two successive flags $\Phi_{i-1}$ and $\Phi_i$ are adjacent 
with $\Phi\cap\Psi\subseteq\Phi_i$ for all $i$.
Finally, we require the {\em diamond condition}, namely that,
whenever  $F\leq G$, with rank$(F)=j-1$ and rank$(G)=j+1$, there 
are exactly two faces $H$ of rank $j$ such that $F\leq H\leq G$. 

The diamond condition implies that given a flag $\Phi$ of $\po$, for each $i \in \{0, \dots, n-1\}$ there exists a unique flag $\Phi^i\in {\mathcal F(\mathcal P)} $ which is $i$-adjacent to $\Phi$.

Given two faces $F$ and $G$ of a polytope ${\mathcal P}$ such that $F \leq G$, the {\em section} $G/F$ of ${\mathcal P}$ is the set of faces $\{H \in \po \ | F \leq H \leq G\}$, with the induced partial order.  If $F_0$ is a vertex, then the section $F_n/F_0$ is called the {\em vertex-figure\/} of $F_0$. Note that every section $G/F$ of a polytope $\po$ is also a polytope and that $\mathrm{rank}(G/F)=\mathrm{rank}(G) - \mathrm{rank}(F) -1$. 

Let $\po$ and $\qo$ be two $n$-polytopes. An {\em isomorphism} from $\po$ to $\qo$ is a bijection $\gamma: \po \to \qo$ such that $\gamma$ and $\gamma^{-1}$ preserve the order.
An {\em anti-isomorphism} $\delta: \po \to \qo$ is a bijection reversing the order, in which case $\po$ and $\qo$ are said to be {\em duals} of each other, and the usual convention is to denote $\qo$ by $\po^*$. (Note that $({\mathcal P}^*)^*\cong {\mathcal P}$.)
An isomorphism from $\po$ onto itself is called an {\em automorphism} of $\po$. The set of all automorphisms of $\po$ forms a group, its automorphism group, denoted by $\Gamma(\po)$. It is not difficult to see that $\Gamma(\po)$ acts freely on $\fl(\po)$, the set of all flags of $\po$.
An anti-isomorphism from $\po$ to itself is called a {\em duality}. When a duality of $\po$ exists, $\po$ is said to be {\em self-dual}. Note that the set of all dualities is not a group, as the product of two dualities is in fact an automorphism. However, the dualities and automorphisms of $\po$ together do form the {\em extended group} of $\po$, denoted by $\bar{\Gamma}(\po)$.

A polytope $\po$ is said to be  {\em regular} if $\Gamma(\po)$ is transitive on the flags of $\po$.

\section{Groups of chiral polytopes}\label{groups}

Intuitively, a chiral polytope is a polytope that has all rotational symmetries but not the reflections. Before discussing chirality, we will thus take a look at what we call the rotational subgroup of a regular polytope.

Let \(\mathcal{P}\) be a regular \(n\)-polytope. Let \(\Phi = \{F_{-1},F_0,...,F_n\}\) be a flag of ${\mathcal P}$. Let $\rho_i\in \Gamma({\mathcal P})$ be the involution sending $\Phi$ to its $i$-adjacent flag. The group \(\Gamma(\mathcal{P}) = \langle \rho_0,...,\rho_{n-1}\rangle\).
As explained in~\cite[Proposition 2B10]{ARP}, the group \(\Gamma(\mathcal{P}) = \langle \rho_0,...,\rho_{n-1}\rangle\) satisfies the {\it intersection property}, namely 
\begin{equation}
\langle \rho_i : i \in I\rangle \cap \langle \rho_j : j \in J\rangle = \langle \rho_k : k \in I\cap J\rangle \textrm{ for all } I,J\subseteq\{0,\ldots,n-1\}
\end{equation}
and, by~\cite[Proposition 2B11]{ARP},  the group \(\Gamma(\mathcal{P}) = \langle \rho_0,...,\rho_{n-1}\rangle\) satisfies the {\it string property}, namely
\begin{equation}
\rho_i\rho_j = \rho_j\rho_i \textrm{ for all } i,j \in\{0,\ldots,n-1\}\textrm{ with }|i-j|>1.
\end{equation}
Let us define
\begin{center}
\(\sigma_j := \rho_{j-1}\rho_j\)
\end{center}
for every \(j = 1,...,n-1\). Then, those \(\sigma_j\) fix every element of \(\Phi\) except for \(F_{j-1}\) and \(F_j\). Moreover, it is easy to see that \(\sigma_j\) cyclically commutes, or in other words ``rotates", the consecutive \(j\)- and \((j-1)\)-faces of the section \(F_{j+1}/F_{j-2}\), which is a polygon. 
The {\em Schl\"afli type} of $\mathcal P$ is the ordered set $\{p_1, \ldots, p_{n-1}\}$ where $p_j := o(\sigma_j)$ for every $j=1, \ldots n-1$.
Now define the subgroup
\begin{center}
\(\Gamma^+(\mathcal{P}) := \langle \sigma_1,...,\sigma_{n-1}\rangle\)
\end{center}
of \(\Gamma(\mathcal{P})\). \(\Gamma^+(\mathcal{P})\) is called the \textit{rotational subgroup} of \(\Gamma(\mathcal{P})\) and is of index at most \(2\) in \(\Gamma(\mathcal{P})\). The group \(\Gamma^+(\mathcal{P})\) contains every word of even length in \(\Gamma(\mathcal{P})\). If \(\Gamma^+(\mathcal{P})\) is of index exactly \(2\), we say that \(\mathcal{P}\) is \textit{directly regular}. If \(\mathcal{P}\) is of Schl\"afli type \(\{p_1,...,p_{n-1}\}\), then the generators of \(\Gamma^+(\mathcal{P})\) satisfy the following relations:
\begin{center}
\(\sigma_j^{p_j} = 1 \)~~ for \(1\le j \le n-1\).

\((\sigma_j\sigma_{j+1}...\sigma_k)^2 = 1\)~~ for \(1\le j <k \le n-1\).
\end{center}
With this intuitive understanding of the behavior the rotational subgroup of a regular polytope, let us properly define a chiral polytope.

An \(n\)-polytpe \(\mathcal{P}\) is called \textit{chiral} if its group \(\Gamma(\mathcal{P})\) has exactly two orbits on the set of flags of \(\mathcal{P}\) such that adjacent flags belong to different orbits.

Let \(\mathcal{P}\) be a chiral polytope and \(\Phi\) be a base flag of \(\mathcal{P}\). By definition of chirality, there are two different types of flags in \(\mathcal{P}\).  A flag \(\Psi\) is called an \textit{even flag} of \(\mathcal{P}\) (with respect to \(\Phi\)) if there exists a sequence \(\Phi=\Phi_0,\Phi_1,...,\Phi_{2k-1},\Phi_{2k} =\Psi\) of adjacent flags. If we do not count the first flag \(\Phi\), there is an even number of flags in this sequence, hence justifying the terminology. A flag that is not even is said to be \textit{odd}. 

Whenever \(\mathcal{P}\) is chiral with base flag \(\Phi = \{F_{-1},F_0, ..., F_n\}\), there exist elements \(\sigma_i \) in \(\Gamma(\mathcal{P})\) such that \(\sigma_i(\Phi) = (\Phi^{i-1})^{i}\) for \(i = 1,...,n-1\) since \(\Phi\) and \((\Phi^{i-1})^{i}\) are in the same orbit. For the sake of simplicity, the flag \((\Phi^{i})^{j}\) will be noted \(\Phi^{i,j}\) from now on. Let us define, for any \(1 \le i \le j\le n-1\),
\begin{equation}
 \tau_{i,j} := \sigma_i\sigma_{i+1}...\sigma_j
\label{tau}
\end{equation}
It can easily be proven that those \(\tau_{i,j}\) are such that \(\tau_{i,j}(\Phi)= \Phi^{i-1,j} \)

Here is a characterization of the chirality, due to Egon Schulte and Asia Weiss.
\begin{thm}
\label{theorem:chiral}\cite{Chiral}
An \(n\)-polytope \(\mathcal{P}\) is chiral if and only if it is not regular, but for some base flag \(\Phi = \{F_{-1},F_0, ..., F_n\}\) of \(\mathcal{P}\) there exist automorphisms \(\sigma_1,...,\sigma_{n-1}\) of \(\mathcal{P}\) with the following properties:
\begin{itemize}
\item \(\sigma_j\) fixes each face in \(\Phi \setminus\{F_{j-1},F_j\}\) for every \(j = 1,.., n-1\).

\item \(\sigma_j\) cyclically permutes consecutive \(j\)-faces in the \(2\)-section \(F_{j+1}/F_{j-2}\) of \(\mathcal{P}\) for every \(j=1,...,n-1\).
\end{itemize}
\end{thm}

Consequently, as soon as we have a chiral polytope \(\mathcal{P}\) and a base flag \(\Phi = \{F_{-1},F_0,...,F_n\}\) of \(\mathcal{P}\), we also have a set \(\{\sigma_1,...,\sigma_{n-1}\}\) of abstract rotations as in Theorem~\ref{theorem:chiral}. Those \(\sigma_i\)'s, however, are not uniquely defined. We can indeed always replace \(\sigma_i\) by its inverse.
We will thus impose an "orientation" on the \(\sigma_i\). Let us denote by \(F_i'\) the only flag such that \(F_{i-1} \le F_i'\le F_{i+1}\) and \(F_i' \ne F_i\). Then, for each \( i :=1,...,n-1\), we will ask that \(\sigma_i(F_i') = F_i\) and thus that, automatically, \(\sigma_i(F_{i-1}) = F_{i-1}'\). The corresponding set  \(\{\sigma_1,...,\sigma_{n-1}\}\) is called the set of \textit{distinguished generators} of \(\mathcal{P}\).
We will now look at some interesting results for chiral polytopes.
\begin{prop}[\cite{Chiral}, Proposition 4]
\label{genchiral}
Let \(\mathcal{P}\) be a chiral polytope. Then \(\Gamma(\mathcal{P}) =\langle \sigma_1,...,\sigma_{n-1}\rangle\). Furthermore, the generators \(\langle \sigma_1,...,\sigma_{n-1}\rangle\) satisfy :
\begin{equation}
\left \{
\begin{array}{ll}
\sigma_i^{p_i} = 1~~ for ~1\le i\le n-1,\\
(\sigma_i\sigma_{i+1}...\sigma_j)^2=1~~ for~ 1 \le i< j \le n-1
\end{array}
\right.
\end{equation}
with the \(p_i\)'s given by the Schl\"afli type of \(\mathcal{P}\).
\end{prop}
\begin{prop}[\cite{Chiral}, Proposition 4]
The faces and cofaces of \(\mathcal{P}\) are directly regular or chiral polytopes.
\end{prop}

Let \(C\) be a chain of a chiral polytope $\mathcal{P}$ and denote by \(\Gamma(C,\mathcal{P})\) the stabilizer of \(C\) in \(\Gamma(\mathcal{P})\). It can then be proved that the stabilizer of the chain \(C:= \{F_i~|~i\notin I\}\) is given by \(\Gamma(C_I,\mathcal{P}) = \langle \tau_{r,s}~|~r\le s, r-1,s \in I\rangle\). We obtain the following proposition as an immediate consequence.
\begin{prop}[\cite{Chiral}, Proposition 6]
The stabilizers of the faces in \(\Phi\) are as follows:
\begin{itemize}
\item \(\Gamma(F_0,\mathcal{P}) = \langle \sigma_2,...,\sigma_{n-1}\rangle\)
\item \(\Gamma(F_i,\mathcal{P}) = \langle \{\sigma_j~|~j\ne i,i+1\} \cup \{\sigma_i\sigma_{i+1}\}\rangle \) for \(i = 1,...,n-2\)
\item \(\Gamma(F_n-1,\mathcal{P}) = \langle \sigma_1,...,\sigma_{n-2}\rangle \)
\end{itemize}
\end{prop}
Similarly to the regular case, groups of chiral polytopes also have an intersection property.
\begin{prop}[\cite{Chiral}, Proposition 7]
Fix \(I,J\subseteq \{-1,0,1,...,n\}\).Then, the following intersection property holds:
\begin{equation}
\begin{split}
\langle \tau_{r,s}~|~r\le & s, r-1,s \in I\rangle  \cap \langle \tau_{r,s}~|~r\le s, r-1,s \in J\rangle\\
& = \langle \tau_{r,s}~|~r\le s, r-1,s \in I\cap J\rangle.
\end{split}
\end{equation}
\end{prop}

\section{Classifying chiral polytopes}\label{algo}

Our goal in this section is to design an algorithm that takes a given group \(G\) as input and that will, as output, give a list of all non degenerate chiral polytopes that have the given group \(G\) as automorphism group. We already discussed briefly how to do so in the previous section and we noticed that an important part for having an effective algorithm is to find sufficient conditions on the generators of \(G\) for the associated coset geometry to always be a chiral polytope. We use the results of~\cite{Chiral} to do so, and then translate those results in the \(C^+\)-groups terminology introduced in~\cite{Hypertopes} to obtain new generators, that allow us to design a faster algorithm.

 Let \(n \ge 3\) and \(p_1,...,p_{n-1}\ge 2\) be integers. Consider a group \(A\) generated by elements \{\(\sigma_1\), \(\sigma_2\), ...,  \(\sigma_{n-1}\)\} satisfying the following conditions :
\begin{equation}
\label{generators}
          \begin{cases}
              \sigma_i^{p_i} = 1 \hspace{3cm}  if~  i = 1, ..., n-1 \\
               (\sigma_i  \sigma_{i+1} ...   \sigma_j)^2 = 1 \hspace{1.5cm} if~ 1  \leq i \leq j \leq n-1 \\
               
            \end{cases}
\end{equation}
where we assume that \(p_i\) is the exact order of \(\sigma_i\) for every \(i\). We note that, since \((\sigma_i\sigma_{i+1})^2 = 1\), we also have that \((\sigma_{i+1}\sigma_i)^2 = (\sigma_i^{-1}(\sigma_i\sigma_{i+1})\sigma_i)^2 = 1\).

We also define
\begin{center}
\(\sigma_0:=\sigma_n := 1\)
\end{center}
and, for \(1  \leq i \leq j \leq n-1\), we put
\begin{equation}
\label{tau_ij}
 \tau_{i,j} := \sigma_i \sigma_{i+1}... \sigma_j .
\end{equation}
Then, \(\tau_{i,i} = \sigma_i\) and \((\tau_{i,j})^2= 1\). We also put \(\tau_{0,j} :=\tau_{i,n} := 1\) .

Lastly, we define 
\begin{center}
\( A_I :=  \langle \tau_{r,s}|  r \leq s ~ and~ r-1, s \in I \rangle\).
\end{center}
Therefore, we have that:
\[
          \begin{cases}
              A_{\emptyset} = A_{\{i\}} = 1\\
              A_{\{\sigma_1,...,\sigma_{n-1}\}} = A\\
              A_{\{\sigma_i,\sigma_{i+1},...,\sigma_j\}} = \langle \sigma_i,\sigma_{i+1}...\sigma_j \rangle
               
            \end{cases}
\]

From now on, we ask that the group A together with its set of generators  \{\(\sigma_1\), ...,  \(\sigma_{n-1}\)\} has the following {\em intersection property} :
\begin{equation}
\label{IC+}
A_I \cap A_J = A_{I \cap J} \hspace{2cm} for~ I,J \subseteq \{ -1, 0, 1, ..., n\}
\end{equation}
We can now state Theorem \(1\) of~\cite{Chiral}, without giving the proof here.
\begin{thm}[\cite{Chiral}, Theorem 1]
\label{schulte}
 Let \(n \ge 3\), \(2\le p_1, ..., p_{n-1} \leq \infty  \) and \( A := \langle \sigma_1, ..., \sigma_{n-1}\rangle\) be a group with properties (\ref{generators}) and (\ref{IC+}). Let \(\mathcal{P} =\mathcal{ P}(A) \) the poset corresponding to \(A\) (the coset geometry associated to \(A\) in our vocabulary). Then \(\mathcal{P}\) has the following properties:
\begin{itemize}

\item \(\mathcal{P}\) is a chiral or directly regular polytope with \(\Gamma^+(\mathcal{P}) = A\) and \(\{\sigma_1,...,\sigma_{n-1}\}\) is the set of distinguished generators of \(\Gamma^+(\mathcal{P})\) associated to some base flag of \(\mathcal{P}\).
\item \(\mathcal{P}\) is of type \(\{p_1,...,p_{n-1}\}\). 
\item \(\mathcal{P}\) is directly regular if and only if there exists an involutory group automorphism \(\rho : A \rightarrow A\) such that \(\rho(\sigma_1) = \sigma_1^{-1},~ \rho(\sigma_2) = \sigma_1^2\sigma_2 ~and~\rho(\sigma_i) = \sigma_i \) for \( i = 3, ..., n-1\)
\end{itemize}
\end{thm}

The polytope \(\mathcal{P}\) is of type \(\{p_1,...,p_{n-1}\}\), according to the theorem. We can, as always, suppose that \(p_i \ne 2\) for otherwise, the polytope we obtain is degenerate, that is, it can be obtained as a direct sum of two smaller polytopes. With that in mind, we can now construct an algorithm that searches for all sets of non involutory elements of the given group \(G\) that generate \(G\) and such that they satisfy properties (\ref{generators}) and (\ref{IC+}). We then have to check that there is no involutory group automorphism \(\rho : A \rightarrow A\) such that \(\rho(\sigma_1) = \sigma_1^{-1},~ \rho(\sigma_2) = \sigma_1^2\sigma_2 \) and \(\rho(\sigma_i) = \sigma_i \) for \( i = 3, ..., n-1\), otherwise the theorem tells us that \(\mathcal{P}\) is not chiral bur regular.
We know that, by doing so, we find all the polytopes we were searching for, since the automorphism group of a chiral polytope is always generated by elements respecting conditions (\ref{generators}) and (\ref{IC+}). This algorithm is the one that was presented in~\cite{HHL}.

 For technical reasons that we will clarify later, this algorithm is quite slow. The theory of \(C^+\)-groups (see~\cite{Hypertopes}) gives us a different characterization of the set of elements generating a chiral polytope. These new generators allow us to build a much faster algorithm. For that reason, we now describe those new set of generating elements and show that there is a bijection from the set of the classical generators to the one coming from \(C^+\)-groups. We will therefore be assured, without having to prove everything again, that the results of Theorem~\ref{schulte} are still valid and that we therefore still have an exhaustive algorithm to find chiral polytopes having a given group \(G\) as automorphism group.

The idea is to take the set of generators of a \(C^{+}\)-group and put additional conditions on them that will force the digon diagram of the geometry associated to be linear. For that, let $G$ be a group with a set of generators \( R := \{ \alpha_1, ..., \alpha_{n-1}\}\), set \( \alpha_0 := 1\) and define 
\begin{center}
\( \alpha_{i,j} := \alpha_i^{-1} \alpha_j \) \hspace{1.5cm} for \( 0 \leq i,j \leq n-1\)
\end{center}
We will also define
\begin{center}
\(G_J^+ := \langle\alpha_{i,j} |  i,j \in J\rangle \)\hspace{1.5cm} for \( J \in \{0, ..., n-1\}\)
\end{center}
And we will naturally ask that this group follows the intersection condition (\(IC^+\)) of the \(C^+\)-groups :
\begin{equation}\label{icplus}
G_J^+\cap G_K^+ = G_{J\cap K}^+ \forall J,K \subseteq \{0, ..., n-1\} \textrm{ with } |J|, |K| \geq 2.
\end{equation}
Up to here, we only asked \((G^+, R)\) to be a \(C^+\)-group. By definition, a {\em $C^+$-group} is a pair $(G^+,R)$ where $G^+$ is a group and $R$ is a set of generators of $G^+$ satisfying (\ref{icplus}). We now have to put more conditions on the generators \(\{\alpha_0,\alpha_1,...,\alpha_{n-1}\}\) so that the diagram of the associated coset geometry is linear. With a bit of intuition, we can remark that the abstract rotation linked to the residue of type \(\{i,j\}\) of the base chamber in the associated coset geometry of (\(G^+,R\)) is \(\alpha_{i,j} = \alpha_i^{-1}\alpha_j\). It is therefore natural to ask that 
\begin{equation}
\label{linear}
\begin{cases}
o(\alpha_i^{-1}\alpha_j) = 2 ~~if~ |i-j| \ge 2 \\
o(\alpha_i^{-1}\alpha_j) \ne 2 ~~if~ |i-j| = 2.
\end{cases}
\end{equation}
This means, in particular, that \(o(\alpha_1) = o(\alpha_0^{-1}\alpha_1) \ne 2\) and that \(o(\alpha_i) =o(\alpha_0^{-1}\alpha_i) = 2\) for all \(i = 2,3,..., n-1\).
We therefore only have to search for a set \(\{\alpha_1,\alpha_2,...,\alpha_{n-1}\}\) in which \(\alpha_1\) is a non involutory element and \(\alpha_i\) is an involution for every \(i=2,3,...,n-1\). As we will see in the section dedicated to the speed analysis of the algorithm, there are usually a lot less involutions than non involutory elements in an almost simple group \(G\). This is the main reason why an algorithm based on these new generators is much faster than an algorithm based on the generators described in~\cite{Chiral}.

Finally, as shown in Theorem~\ref{chirality} below, we need to check that there is no element inverting all generators for otherwise the polytope is regular instead of chiral.

\begin{thm}[\cite{Hypertopes}, Theorem 8.2]
\label{chirality}
Let \((G^+,R)\) be a \(C^+\)-group. Let \(\Gamma = \Gamma(G^+,R)\) be the coset geometry associated to (\(G^+,R\)). If \(\Gamma\) is a hypertope and \(G^+\) has two orbits on the set of chambers of \(\Gamma\), then \(\Gamma\) is chiral if and only if there is no element of \(G^+\) that inverts every generators of \(R\). On the other hand, if such an element \(\sigma \in Aut_I(G^+)\), that inverts every generators of \(R\), exists, then the group \(G^+\) extended by \(\sigma\) is regular on \(\Gamma\).
\end{thm}

\subsection{Equivalence}
\label{Equivalence}
Let $G$ be a group. Define \(A \) to be the set of all sets \(R :=\{\sigma_1,...,\sigma_{n-1}\}\) such that \(G = \langle \sigma_1,...,\sigma_{n-1} \rangle\) and such that \(R\) satisfies conditions~(\ref{generators}) and~(\ref{IC+}) and define \(B\) to be the set of all sets \(R' := \{\alpha_1,...,\alpha_{n-1}\}\) such that \((G,R')\) satisfies conditions~(\ref{icplus}) and~(\ref{linear}). We also suppose that every \(\sigma_i\) appearing in \(A\) has an order strictly greater than two. This simply asks that the diagram associated is connected.

We now show that there is a bijection between the two sets \(A\) and \(B\).
\begin{thm}
There exist a bijection \(\varphi : A\rightarrow B\) given by \[\varphi(\{\sigma_1,...,\sigma_{n-1}\}) =\{\sigma_1,\sigma_1\sigma_2, ..., \sigma_1\sigma_2...\sigma_{n-1}\}.\]
\end{thm}
\begin{proof}Let us define \[\psi: B \to A: \psi(\{\alpha_1,...,\alpha_{n-1}\})=\{\alpha_1,\alpha_1^{-1}\alpha_2, \alpha_2^{-1}\alpha_3, ...,\alpha_{n-2}^{-1}\alpha_{n-1}\}. \]
We first prove that \(\varphi\) and \(\psi\) are well defined, meaning here that their images are respectively in \(B\) and \(A\) and then we will see that \(\varphi\circ \psi = Id _B\) and \(\psi \circ\varphi = Id_A\).

Let us check that \(\{\sigma_1,\sigma_1\sigma_2, ..., \sigma_1\sigma_2...\sigma_{n-1}\}\) is in \( B\) for any \(\{\sigma_1,...,\sigma_{n-1}\}\) in \( A\).

\begin{itemize}
\item The order of \(\sigma_1\) is strictly bigger than $2$, by definition, as noted in the introduction of this section. The order of \(\sigma_i...\sigma_j\) is equal to $2$ by the property~(\ref{generators}) of \(A\). The generators are thus all involutions except the first one, as it should be. Similarly, we also have that \((\sigma_1...\sigma_i)^{-1} \sigma_1...\sigma_i...\sigma_j = \sigma_{i+1}...\sigma_j\) so that all the relations of~(\ref{linear}) are automatically satisfied by the condition~(\ref{generators}).
\item Since \( G = \langle\sigma_1, ..., \sigma_{n-1}\rangle \) it is also generated by \( \{\sigma_1, \sigma_1\sigma_2, ..., \sigma_1\sigma_2...\sigma_{n-1}\}\).
Indeed, we have \(\sigma_2 =\sigma_1^{p_1-1}(\sigma_1\sigma_2)\) and similarly we can get \(\sigma_3, ..., \sigma_{n-1}\).
\item We now have to check that~(\ref{icplus}) holds. To this end, we need to consider all the groups generated by some of the \(\alpha_{i,j} = \alpha_i^{-1}\alpha_j\) where the \(\alpha_i\)'s are the distinguished generators of \(G^+\). To avoid confusion in the notation,  let us define \( \sigma_{i,j} :=( \sigma_1...\sigma_i)^{-1}(\sigma_1...\sigma_j)\) for \(1\le i \le j \le n-1\) (i.e: the \(\sigma_{i,j}\)'s are just the \(\alpha_{i,j}\)'s expressed with the generators given by \(\varphi(\{\sigma_1,...,\sigma_{n-1}\})\).Then 
\begin{center}
\(\sigma_{i,j} = \left\{
    \begin{array}{ll}
        1 &\mbox{if \(i=j\)} \\
        \sigma_i & \mbox{if \( i+1 = j\)}\\
          \sigma_i...\sigma_j & \mbox{if  \(i+2 \leq j\)}.
    \end{array}
\right.
\)\\
\end{center}

We can now see that \(\sigma_{i,j} = \tau_{i+1,j} \) where \(\tau_{i,j}\) is defined as in~(\ref{tau_ij}) with the additional convention that \(\tau_{i+1,j} = 1 \) if \(i+1> j\). We are thus ready to check~(\ref{icplus}). Indeed, for any \( J \in \{0,1,..., n-1\}\), we have \( G_J^+ := \langle \sigma_{i,j}| i,j \in J\rangle  = \langle \sigma_{i,j}|i \leq j, i,j \in J\rangle =  \langle \tau_{i+1,j}|i \leq j, i,j \in J \rangle = \langle \tau_{i+1,j}|i+1\leq j, i,j\in J\rangle = \langle \tau_{k,j}|k\leq j, k-1,j \in J\rangle = A_J\) by definition of \(A_J\).\\
Hence, for all \(J,K \in \{0,1,..., n-1\}\), we have that \(G_J^+\cap G_K^+ = A_J\cap A_K = A_{J\cap K}= G_{J\cap K}^+\). 
\end{itemize}

Let us now check that \(\psi(\{\alpha_1,...,\alpha_{n-1}\}) = \{\alpha_1,\alpha_1^{-1}\alpha_2, \alpha_2^{-1}\alpha_3, ...,\alpha_{n-2}^{-1}\alpha_{n-1}\} \) is in A for every \(\{\alpha_1,...,\alpha_{n-1}\}\in B\). 
\begin{itemize}
\item Once again we first need to verify that condition~(\ref{generators}) holds. 
Since \[\alpha_i^{-1}\alpha_{i+1}\alpha_{i+1}^{-1}\alpha_{i+2}...\alpha_{j-1}^{-1}\alpha_j = \alpha_i^{-1}\alpha_j,\] condition~(\ref{linear}) implies that the generators satisfy every relation of (\ref{generators}).

\item \( G = \langle\alpha_1,...,\alpha_{n-1}\rangle \) by definition of \(B\). Thus, since \(\alpha_i = \alpha_1(\alpha_1^{-1}\alpha_2)...\) \((\alpha_{i-1}^{-1}\alpha_i)\), it comes that \( G= <\alpha_1,\alpha_1^{-1}\alpha_2, \alpha_2^{-1}\alpha_3, ...,\alpha_{n-2}^{-1}\alpha_{n-1}>\).

\item Similarly as in the previous case, we start by taking a look at the expression of the \(\tau_{i,j}\) derived by the given generators. Here \[\tau_{i,j} = (\alpha_{i-1}^{-1}\alpha_i)(\alpha_i^{-1}\alpha_{i+1})...(\alpha_{j-1}^{-1}\alpha_{j}) = \alpha_{i-1}^{-1}\alpha_j = \alpha_{i-1}^{-1}\alpha_j=: \alpha_{i-1,j}.\]
Therefore, for every \( J \in \{0,1,...,n-1\}\)  we have 
\begin{align*}
A_J = \langle \tau_{i,j}~|~i \leq j, i-1,j \in J \rangle &=\langle \alpha_{i-1,j}~|~i \leq j ,i-1,j \in J \rangle \\  &=\langle\alpha_{k,j}~|~k+1 \leq j, k,j \in J\rangle \\&= \langle \{\alpha_{k,j}~|~ k+1\leq j, k,j\in J\} \cup \{\alpha_{j,j} = 1\}\rangle \\&= \langle\alpha_{k,j}~|~k \leq j, k,j \in J\rangle \\&= G_J^+.
\end{align*}
Hence, for all \(  J,K \subseteq \{0,1,..., n-1\}\), we have \(= A_J\cap A_K =G_J^+\cap G_K^+ = G_{J\cap K}^+ = A_{J\cap K}\). 

Technically, in the condition~(\ref{IC+}), we also have to check the property for the subsets \(J, K \subseteq \{-1,0,1,...,n\}\) containing \(-1\) and \(n\).
For that, it suffices to notice that, since \(\tau_{0,i} =\tau_{i,n} = 1\) for all \( i \), we have, for all \(\ J \in  \{0,1,..., n\}\), that \( A_j =A_{J\cup \{-1\}} = A_{J\cup \{n\}} = A_{J\cup \{-1,n\}}\).
\end{itemize}
Finally, it remains to verify that \(\varphi\) and \(\psi\) are inverse. This is straightforward from their respective definitions.
\end{proof}

We now have all the tools to write our algorithm. The only detail we did not take into account yet is to verify that the criterion for the polytope to be indeed chiral, and not regular, is the same in both cases. 

\begin{lemma}
Let \(G\), \(A\), \(B\), \(\varphi\) and \(\psi\) be as before. Then, for every \(\{\sigma_1,...,\sigma_{n-1}\}\) \( \in A\) and for every \(\{\alpha_1,...,\alpha_{n-1}\}\in B\), we have that:

There exists an involutory element \(\rho \in G\) such that \(\rho(\sigma_1)= \sigma_1^{-1}\) , \(\rho(\sigma_2) = \sigma_1^2\sigma_2\) and \(\rho(\sigma_i) = \sigma_i \) for every \( i = 3,...,n-1\) if and only if there exists a \(g\in G\) such that \(g\) inverts every elements of \(\varphi(\{\sigma_1,...,\sigma_{n-1}\})\).
\end{lemma}
\begin{proof}
Actually, the elements \(\rho\) and \(g\) are the same.

Let us first suppose that we have \(\rho\) such that \(\rho(\sigma_1)= \sigma_1^{-1}\) , \(\rho(\sigma_2) = \sigma_1^2\sigma_2\) and \(\rho(\sigma_i) = \sigma_i \) for every \( i = 3,...,n-1\). Then, this \(\rho\) is also such that \(\rho(\sigma_1)= \sigma_1^{-1}\) and \(\rho(\sigma_1\sigma_2) = \rho(\sigma_1)\rho(\sigma_2) = \sigma_1^{-1} \sigma_1^2\sigma_2 = \sigma_1 \sigma_2\). Similarly, we have that \(\rho(\sigma_1...\sigma_j) =\sigma_1...\sigma_j\). This is exactly what we wanted since every \(\sigma_1...\sigma_j\) is an involution.

Conversely, if we have \(g\in G\) such that it inverts every elements of \(\varphi(\{\sigma_1,...,\sigma_{n-1}\})\), it means that \(g(\sigma_1) = \sigma_1^{-1}\) and \(g(\sigma_1...\sigma_i) = \sigma_1...\sigma_i\). In particular, it also means that \(g(\sigma_1\sigma_2) = g(\sigma_1)g(\sigma_2) = \sigma_1^{-1}g(\sigma_2)\) so that \(g(\sigma_2) = \sigma_1^{2}\sigma_2\). Similarly \(g(\sigma_3) = g(\sigma_2)^{-1}g(\sigma_1)^{-1}g(\sigma_1\sigma_2\sigma_3) = \sigma_2^{-1}(\sigma_1^{-2}\sigma_1)(\sigma_1\sigma_2\sigma_3) = \sigma_3\) and so on for the others.
\end{proof}

\subsection{Code}

In this section, we present in more details the code that has been written. We have proved in the previous section that the algorithm that is searching for all generators of a \(C^+\)-group \(G\) with the additional condition on the linearity of the diagram given by~(\ref{linear}) is correct in the sense that it will find all existing chiral polytopes having \(G\) has an automorphism group. 

Let \(G\) be a group and let \(g,g' \in G\). Recall that the elements \(g\) and \(g'\) are said to be {\em conjugate} if there exists an element \(h\in G\) such that \( g' = hgh^{-1}\). The relation \(g\sim g'\) if and only if \(g\) and \(g'\) are conjugate is an equivalence relation and therefore partitions \(G\) into equivalence classes. We will denote the conjugacy class of \(g\in G\) by \(Cl(g) := \{g' \mid g\) and \(g'\) are conjugate\(\}\).
Let \(\alpha, \beta\) and \(g\) be elements of \(G\). Then if \(\beta = g\alpha g^{-1}\), we have that \(o(\alpha) = o(\beta)\).

\begin{lemma}
Let \(B\) be the set defined in Section~\ref{Equivalence}.
Let  \(R = \{\alpha_1,...,\alpha_{n-1}\}\) be a set of elements of \(G\) and \(g\) be an element of \(G\). Then,  \(R = \{\alpha_1,...,\alpha_{n-1}\}\) is in \(B\) if and only if \(R' = \{g\alpha_1 g^{-1},g\alpha_2 g^{-1},..., g\alpha_{n-1}g^{-1}\}\) is in \(B\) and the chiral polytopes associated to both sets of generators are isomorphic.
\end{lemma}
\begin{proof}
The first part of the proof directly comes from the definition of \(B\) and the above remark on the equality of the orders of two conjugated elements.

 Let \(\Gamma(G,R) = \Gamma(G,(G_i^+)_{i\in I})\) and \(\Gamma(G,R') =\Gamma(G,(H_i^{+})_{i \in I})\) be the coset geometries associated to \((G,R)\) and \((G,R')\) using~\cite[Construction 8.1]{Hypertopes}.  
Our goal is to prove that \(G_i^+\) and \(H_i^+\) are conjugated for every \(i\in I\), so that the two coset geometries are isomorphic.

Let us put \(\alpha_0 := 1_G = g1_Gg^{-1}\) and define \(\alpha_{i,j} = \alpha_i^{-1}\alpha_j\) and \(\alpha_{i,j}' = (g\alpha_ig^{-1})^{-1}(g\alpha_jg^{-1}) = g \alpha_i^{-1}\alpha_j g^{-1} =g\alpha_{i,j}g^{-1}\) for every \(i,j \in \{0,1,...,n-1\}\). We thus notice that \(\alpha_{i,j}\) and \(\alpha_{i,j}'\) are conjugate. It is now sufficient to notice that \(G_k^+ = \langle \alpha_{i,j}~|~i,j \ne k\rangle\) and \(H_k^+ = \langle \alpha_{i,j}'~|~i,j \ne k\rangle\). We therefore see that \(G_k^+ = gH_k^+g^{-1}\) for every \(k = 0,1,...,n-1\).
By the above remark, this concludes the proof.
\end{proof}

This lemma shows that we can restrict our choice of \(\alpha_1\) to one element of each conjugacy classes of elements of order greater than or equal to 3, since every chiral polytope obtained by taking a second element \(\tilde{\alpha_1}\) that lies in the same conjugacy class as \(\alpha_1\) is isomorphic to a chiral polytope already listed thanks to \(\alpha_1\).

Here under is the main section of the {\sc Magma}~\cite{Magma} code we developed. The details about the different functions can be found in the following paragraph. We begin by selecting one element in each conjugacy classes of order strictly greater than two and listing them in a variable called \(K\). We also store in $I$ all involutions in the group. Then, for each element \(x\) of \(K\), we call the function ``OrderedSubsetsGroups()" which will give as an output, every possible combination of fitting generators with \(x\) as the first element. We then group all the outputs in the variable "gen" and we verify, one by one, that the associated geometry is chiral and not isomorphic to one of the chiral polytopes already found previously.

\begin{verbatim}
function ChiralPolytopes(G) 
    c:= ConjugacyClasses(G);
    K:=[[c[i][3]] : i in [1..#c] | c[i][1] gt 2 and c[i][3] in G];
    I := [];
    i := 2;
    while c[i][1] eq 2 do
      I cat:= [c[i][3]^x : x in Transversal(G,Centralizer(G,c[i][3]))];
      i := i + 1; 
    end while;
    Inv:= [i : i in [1..#I]];
    gen := [];
    for i:= 1 to #K do
        subsets:=[];
        s := K[i][1];
        S := [];
        OrderedSubsetsGroups(s,~S,~Inv,~subsets,~G,~K,~I);
        gen cat:=[subsets];
    end for;
    Chiral_poly := [];
    for i := 1 to # gen do
        for j := 1 to #gen[i] do 
            new := true;
            x := [];
            for k := 1 to #gen[i][j] do
                Append(~x,I[gen[i][j][k]]);
            end for;
            poly := K[i] cat x;
            H := sub<G|poly>; 
            if #H eq #G then 
                polybis := [poly[i]^-1: i in [1..#poly]];
                // Here we check if we get a chiral polytope or not
                if not(IsHomomorphism(H,H,polybis)) then
                    for known in Chiral_poly do
                        N := sub<G|known>;
                        if #known eq #poly then
                            if #H eq #sub<G|known> and
                                IsHomomorphism(H,N,known) then
                                    new := false;
                                    break known;
                            end if;
                        end if;
                    end for;
                    if new then
                        Append(~Chiral_poly,poly);
                    end if;
                end if;
            end if;
        end for;
    end for;
    for i := 1 to #Chiral_poly do
        DiagramPlus(G,Chiral_poly[i]);
    end for;	
    return Chiral_poly;
end function;
\end{verbatim}

Below are the details about the recursive function finding all suitable set of generators having a fixed element \(x\) as first generator. Notice that we are checking the properties that we mentioned earlier, for example the intersection property and conditions on the orders of the \(\alpha_i^{-1}\alpha_j\)'s. The restriction to the normalizers and centralizer is a straightforward optimization deduced once again from property~(\ref{linear}) of the generators.

\begin{verbatim}
procedure OrderedSubsetsGroups (s,~S,~Inv,~subsets,~G,~K,~I)
    U := RingOfIntegers();
    if #Inv gt 0 then
        x := [I[S[i]]: i in [1..#S]];
        if #sub<G|s,x> ne #G then
            if #S eq 1 then
                Sub_s := sub<G|s>;
                Norm := Normalizer(G,Sub_s);
                Inv:= [i: i in Inv|I[i] in Norm];
            end if;
            if #S ge 2 then
                Centr := Centralizer(G,I[S[#S-1]]);
                Inv := [i: i in Inv| I[i] in Centr];
            end if;
            if #S ge 1 then
                Inv2 := Inv;
                Centr := Centralizer(G,I[S[#S]]);
                Inv := [i : i in Inv| not I[i] in Centr];
            end if;
            for i:= 1 to #Inv do
                y := true;
                if #S ge 1 then
                    if Order(s^-1*I[Inv[i]]) ne 2 then
                        y := false;
                    end if;
                end if;
                if y then
                    M := sub<G|s,x,I[Inv[i]]>;
                    if not HasIntersectionPropertyPlus(M) then
                        y := false;
                    end if;
                end if ;
                if y then
                    if #S eq 0 then
                        X := Append(S,Inv[i]);
                        Y := Exclude(Inv, Inv[i]);
                        OrderedSubsetsGroups(s,~X,~Y,~subsets,~G,~K,~I);
                    end if;
                    if #S ge 1 then
                        X := Append(S,Inv[i]);
                        Y := Exclude(Inv2, Inv[i]);
                        OrderedSubsetsGroups(s,~X,~Y,~subsets,~G,~K,~I);
                    end if;
                end if;
            end for;
        end if;
    end if;
    Append(~subsets,S);	
end procedure;
\end{verbatim}
Below is a function that allows the output to be easily readable and understood. It mainly shows in the output the Schl\"afli type of the polytope given as input.
\begin{verbatim}
function DiagramPlus(G,gens)
    schlafli := [];
    Append(~schlafli, Order(gens[1]));
    for i := 2 to #gens do
        Append(~schlafli, Order(gens[i-1]^-1*gens[i]));
    end for;
    print "New chiral of type ", [y : y in schlafli],"for group of order ",#G;
    return schlafli;
end function;
\end{verbatim}

\subsection{Speed comparison and complexity analysis}

An algorithm based on the generators described in~\cite{Chiral} already exists. Please refer to \cite{HHL} for an implementation of this algorithm. Our goal in creating this new algorithm was to implement a faster algorithm, so that we could find results for groups that were too big for the old algorithm. The new algorithm is indeed considerably faster than its predecessor when dealing with almost simple groups. 

Suppose that \(G\) is an almost simple group and denote by \(i_2(G) := \{g\in G \mid o(g) = 2\}\) the subset containing all the involution of \(G\). Our goal is to quantify the ratio \(i_2(g) / |G|\) and to get a bound for it. Let us state Corollary \(4.5\) of~\cite{Liebeck} without proof.

\begin{prop}[\cite{Liebeck}, Corollary 4.5]
Let \(G\) be an almost simple group. Then there is a constant c such that \(i_2(G) < c|G|^{\frac{5}{8}}\) unless \(soc(G)= L_2(q)\).
\end{prop}

Therefore, if \(soc(G) \ne L_2(q)\), we have that \(i_2(g) / |G| < \frac{ c|G|^{\frac{5}{8}}}{|G|} = c |G|^{-\frac{3}{8}}\) for some constant \(c\). This ratio is thus always really small and tends to \(0\) when the cardinality of \(G\) goes to infinity.

Table~\ref{timings} gives a speed comparison of both algorithms. The left column contains the almost simple group \(G\) given as an input to the algorithms, column BLT contains the execution times of our algorithm (BLT stands for Buekenhout-Leemans-Tranchida), column HHL contains the execution times of the Hartley-Hubard-Leemans algorithm. Both algorithms were executed on the same machine having 48 cores running at 2.6Ghz and 384 gigabytes of memory. Observe that {\sc Magma} does not do parallel computing and that most groups required very few gigabytes of ram to be fully analysed with our programs. Column $i_2(G)/\# G$ gives the ratio (number of involutions)/(order of the group) and the last column gives the order of $G$. When question marks appear in the table, it means we decided to stop the computations after more than 18 days of computing. One of the most spectacular timing differences is for the group $PSU(3,7)$ which took less than 20 minutes to be analysed with the BLT algorithm and more than 18 days with the HHL algorithm.

\begin{table}[h]
\begin{tabular}{|c|c|c|l|l|}
\hline
Group&BLT algo&HHL algo&$i_2(G)/\# G$&$\#G$\\
\hline
PSL(3,2)&0.05s&0.2s&0.125&168\\
PSL(3,3)&3s&20s&0.02&5616\\
PSL(3,4)&13s&210s&0.015&20160\\
PSL(3,5)&107s&12867s&0.002&372000\\
PSL(3,7)&1091s&433480s&0.0014&1876896
\\
PSL(3,8)&1990s&????&0.00028&16482816\\
PSL(3,9)&52226s&????&0.00017&42456960\\
PSL(3,11)&49608s&????&0.000076&212427600\\
PSL(3,13)&125831s&????&0.00011&270178272\\
PSL(3,16)&276124s&????&0.000049&1425715200\\
PSL(3,17)&1770068s&????&0.00000000029&6950204928\\

\hline
PSU(3,2)&0.14s&0.06s&0.053&72\\
PSU(3,3)&0.93s&12s&0.011&6048\\
PSU(3,4)&8s&287s&0.0097&62400\\
PSU(3,5)&217s&4155s&0.0014&126000\\
PSU(3,7)&1180s&1617730s&0.0011&5663616\\
PSU(3,8)&1512s&????&0.00065&5515776\\
PSU(3,9)&12144s&????&0.00014&42573600\\
\hline
Alt(5)&0.03s&0.14s&0.25&60\\
Alt(6)&1.8s&6s&0.125&360\\
Alt(7)&16s&45s&0.042&2520\\
Alt(8)&403s&1222s&0.017&20160\\
Alt(9)&15879s&62917s&0.0073&181440\\
Alt(10)&480012s&????&0.0030&1814400\\
\hline
Sym(5)&0.17s&0.72s&0.21&120\\
Sym(6)&5s&18s&0.10&720\\
Sym(7)&45s&159s&0.046&5040\\
Sym(8)&1319s&7560s&0.019&40320\\
Sym(9)&49325s&891405s&0.0072&362880\\
Sym(10)&1751937s&????&0.0026&3628800\\
\hline
\end{tabular}
\centering 
\caption{Time comparison between both algorithms for some key groups.\label{timings}}
\end{table}

The main reason why the new algorithm is much faster is, as already discussed earlier, that it relies almost entirely on finding suitable involutions, instead of having to go through every non involutive elements in the group. The new algorithm is therefore significantly better as soon as the given group \(G\) contains few involutions compared to non involutory elements. 

As a final remark, we observe that the BLT algorithm actually computes all the polytopes (regular or chiral) that have the given group $G$ as rotational subgroup. In the case where $G$ is a simple group, every abstract regular polytope on which $G$ act regularly is such that $G$ is also the rotational subgroup. Hence the BLT algorithm not only computes all chiral polytopes but also all regular polytopes for $G$. In order to achieve that with the code given in this paper, it suffices to remove the test for chirality in the function {\sf ChiralPolytopes} above.

\end{document}